\documentclass[12pt]{article} 
\usepackage{verbatim,amsmath,amssymb}
\usepackage{epsfig}
\usepackage{graphicx}
\usepackage{index}
\usepackage{amscd}
\input math111.def
\numberwithin{equation}{section}
\newtheorem{theor}{Theorem}[section]
\newtheorem{defi}[theor]{Definition}
\newtheorem{lemma}[theor]{Lemma}

\newtheorem{prop}[theor]{Proposition}

\def\bs{\vskip 0.4cm}
\def\supp{\mbox{\rm supp }}

\newif\ifmarglab
\makeatletter
\def
\label#1{\@bsphack\ifmarglab\marginpar{LAB:#1}\fi\if@filesw {\let\thepage\relax
   \def\protect{\noexpand\noexpand\noexpand}%
   \edef\@tempa{\write\@auxout{\string
      \newlabel{#1}{{\@currentlabel}{\thepage}}}}%
   \expandafter}\@tempa
   \if@nobreak \ifvmode\nobreak\fi\fi\fi\@esphack}
\makeatother

\begin{document}
\title{ Interpolatory Estimates, Riesz Transforms and
  Wavelet Projections}
\author{ Paul F. X . M\"uller \thanks{P.F.X.M.  is supported
by the Austrian Science foundation (FWF) Pr.Nrs. P23987, P22549  and by the Hausdorff Research Institute for Mathematics, Bonn.} and Stefan M\"uller}
\date{August 21, 2014.}
\maketitle

\begin{abstract}
We prove that  directional wavelet projections and Riesz transforms 
are related by interpolatory estimates. The exponents of interpolation
 depend  on the H\"older estimates of the wavelet
system. This paper complements and continues  previous work \cite{smue2,  LMM} on Haar
projections. 

 \end{abstract}


\section{Introduction}

This paper is  concerned with  wavelet systems, 
directional wavelet projections and their 
estimates  in terms of  Riesz
transforms.
We continue and extend the methods introduced  in  \cite{smue2} and
\cite{LMM}.

Let $\cF$ denote the $L^2 ( \bR^n) $ normalized Fourier transform. 
The  Riesz transform 
$R_i$ is the  Fourier multiplier defined by    
\begin{equation}
\cF(R_i(u))(\xi) = -\sqrt{-1} \frac{\xi_i}{|\xi|}\cF(u)(\xi)
\quad\text{where}\quad 1 \le i\le n , \quad \xi = ( \xi_1,\dots , \xi_n ).
\end{equation}

Let $\cS$ denote the collection of dyadic cubes in $ \bR^n ,$
and $\cA = \{0,1\}^n \sm \{(0,0,\dots ,0\} . $
We let 
$$\{ \vp _Q^{(\varepsilon)} : Q \in \cS, \varepsilon \in \cA \} $$
 denote  an admissible  wavelet system of  H\"older exponent 
  $ 0 < \a \le 1  $ and decay estimates of order $\d>0.$ 
(The  definition is  given in \eqref{structure011} below.)
For a fixed direction    $\varepsilon\in \cA $ the  
associated orthogonal  wavelet  projection is defined as 
$$W^{(\varepsilon)} (u) = \sum_{ Q\in \cS} 
\la u ,  \vp_Q^{(\varepsilon)}\ra  \vp_Q^{(\varepsilon)} |Q|^{-1}, \quad\quad u \in  L^2(\bR^n). $$
The results of this paper give pointwise estimates for the directional wavelet
projection  $W^{(\varepsilon)}$ in
terms of the Riesz transforms. 

\subsubsection*{Admissible wavelet systems.}
We specify  now  the wavelet systems  we
use
in this paper. Recall that $I \sb \bR $ is a dyadic interval if there exist 
natural numbers $k, m \in \bZ $ so that $I  = [(k-1)2^m , k 2^m [ . $
Let $I_1 , \dots, I_n $ be dyadic intervals in $\bR $
so that $ |I_i | = |I_j| .$  
Define the dyadic cube $Q \sb \bR^n , $ as 

$$ Q = I_1 \times \cdots \times I_n . $$

We let 
 $s(Q) $ denote the  side length of $Q , $ thus  $s(Q) = |I_1| $.
Let $\cS $ denote the collection of all dyadic cubes in  $\bR^n  $
and put $ \cA = \{ \varepsilon \in \{ 0, 1 \}^n : 
\varepsilon \ne ( 0,\dots . 0 ) \} . $

We say that
$$\{ \vp_Q^{(\varepsilon)} : Q \in \cS, \varepsilon \in \cA \} $$
 is an admissible  wavelet system  
if 
$\{ \vp_Q^{(\varepsilon)} /\sqrt{| Q|} : Q \in \cS, \varepsilon \in \cA \} $
is  an orthonormal basis in $L^2(\bR^n)$
 and
there exists $C > 0 , \d > 0  $ and $ 0 < \a \le 1 $ 
so that the  following  conditions hold:

\begin{subequations} \label{structure011}
\begin{enumerate}
\item Localization with decay estimates:
\begin{equation}\label{structure011a}
 |\vp^{(\varepsilon)}_Q (x)| \le C \left( 1 + \frac{\dist(x, Q) }{s(Q) }
 \right)^{-n(1+\d)},   \quad\quad x \in \bR^n.
\end{equation}
\item H\"older estimates of order $\a :$
\begin{equation} \label{structure011c} 
| \vp^{(\varepsilon)}_Q (x) - \vp^{(\varepsilon)}_Q (t)| 
\le  C s(Q)^{-\a}|x-t|^\a \left( 1 + \frac{\dist(x, Q) }{s(Q) }
\right)^{-n(1+\d)} , 
\quad\quad |x-t | \le  s(Q) .
\end{equation} 
\item Sectional oscillation for  $ i\in \{j \le n : \varepsilon _j = 1 \} :$
\begin{equation}\label{structure011b}
 |\bE_i (\vp^{(\varepsilon)}_Q )(x) | \le C s(Q) \left( 1 + \frac{\dist(x, Q) }{s(Q) }
 \right)^{-n(1+\d)},   \end{equation}
where $ \bE_i $ denotes integration with respect to the variable $ x_i
, $
\begin{equation}\label{int121} 
\bE_{{i} } (f)(x) = \int^{x_{i}}_{-\infty}f(x_1,\dots,s,\dots,  x_n)ds ,
\quad\quad x = ( x_1, \dots,x_i,\dots  x_n).
\end{equation}
\end{enumerate}
\end{subequations}
We refer to $\d > 0 $ and $ 0< \a \le 1 $ as the decay and H\"older
exponents of a  wavelet system satisfying \eqref{structure011}.
\subsubsection*{Directional wavelet projections.}

We fix 
an admissible  wavelet system $\{ \vp_Q^{(\varepsilon)} : Q \in \cS, \varepsilon \in \cA \} $. 
For  a given direction   $\varepsilon\in \cA,$  let  $W^{(\varepsilon)}$
denote the associated  projection on $ L^2(\bR^n), $ 
$$W^{(\varepsilon)} (u) = \sum_{ Q\in \cS} 
\la u ,  \vp_Q^{(\varepsilon)}\ra  \vp_Q^{(\varepsilon)} |Q|^{-1}, \quad\quad u \in  L^2(\bR^n). $$
We summarize next the   main estimates  in \cite{smue2, LMM}, and
relate them to the results of the present paper.

\paragraph{Review of  \cite{smue2}.} 

 If the  H\"older exponent of the   wavelet system satisfies  $0 < \a <1, $ 
the following  Hilbertian estimates
for  $W^{(\varepsilon)}$ are  obtained with the method introduced in
\cite{smue2},
 \begin{equation}\label{19juni14-7}
\|W^{(\varepsilon)}(u)\|_2 \le 
A\left(  \|u\|_2^{1-\a}\|R_ {i_0}(u)\|_2^{\a} + 
 \frac{\|u\|_2^{1-\a} \|R_{i_0} (u)\|_2^\a    - \|R_{i_0} (u)\|_2 }{ 2^{1 -\a }   - 1}   \right),
\end{equation}
whenever $ \varepsilon = (\varepsilon _1, \dots \varepsilon_n ) $ and
$\varepsilon_{i_0} =1 . $ We have 
$ A = A(  \a,\d ) \to \infty$ as $ \a \to 0  $ or $ \d \to 0 .$
The Lipschitz case 
when  $\a = 1 $ is of particular interest.  It appears as the  limit as  $\a \to 1 $ of 
the estimates \eqref{19juni14-7}. By   l'H\^{o}pital's rule  \eqref{19juni14-7} implies
 \begin{equation}\label{4nov011-1}
\|W^{(\varepsilon)} (u)\|_2 \le A(1, \d) \left( 1 + \log \frac{\| u \|_2}
{\|R_{i_0} u\|_2 } \right) 
\|R_{i_0} u \|_2. 
\end{equation}
If    $0 < \a <1, $ is fixed and if one is not interested in the limiting behvior 
as $ \a \to 1 , $ then a simplified form of  \eqref{19juni14-7} is as follows,
\begin{equation}\label{4nov011-2} 
\|W^{(\varepsilon)} (u)\|_2 \le \frac{A ( \a, \d )}{1-\a}  
\| u \|_2^{1-\a}
\|R_{i_0} u \|_2^{\a} ,
\end{equation}
The estimates   \eqref{19juni14-7},  \eqref{4nov011-1} and \eqref{4nov011-2} were proven  in 
\cite{smue2} by  Cotlarization
of the operator  $W^{(\varepsilon)}.$

\paragraph{The  present paper}  extends 
the $L^2(\bR^n) $ estimates \eqref{19juni14-7}  to the scale of $ L^p (\bR^n) $ spaces.
We use below  the abbreviation  $\|R_{i_0}\|_p = \| R_{i_0} : L^p( \bR^n) \to 
L^p( \bR^n) \| .$  Our main result asserts   that  for  H\"older exponents   
$0 < \a <1  $ and  $ 1 < p < \infty, $ 
 \begin{equation}\label{19juni14-8} 
\|W^{(\varepsilon)}(u)\|_p \le 
C \left( 
\|R_{i_0}\|_p^{-\a} \|u\|_p^{1-\a}\|R_ {i_0}(u)\|_p^{\a} + 
 \frac{\|R_{i_0}\|_p ^{1-\a} \|u\|_p^{1-\a} \|R_{i_0} (u)\|_p^\a    - \|R_{i_0} (u)\|_p }{ 2^{1 -\a }   - 1}   \right)
\end{equation}
whenever $ \varepsilon = (\varepsilon _1, \dots \varepsilon_n ) $ and
$\varepsilon_{i_0} =1 . $ 

The asymptotic behavior of the constants
$ C = C(p ,  \a,\d )$
 is as follows,
 $$C(p ,  \a,\d ) = \frac{p^2 C(  \a,\d )}{p-1} \quad\rm{ and }\quad
 C(  \a,\d ) \to \infty \quad \rm{as} \quad  \a \to 0  \quad \rm{or}\quad  \d \to 0 .$$
Again the estimates for  the Lipschitz case  $\a = 1 $ 
appear as limit  of \eqref{19juni14-8} by using   l'H\^{o}pital's rule,
\begin{equation}\label{4nov011-4} 
 \|W^{(\varepsilon)} (u)\|_p \le A( p,\d)  \left( 1 + \log \frac{\| u \|_p\|R_{i_0} \|_p}
{\|R_{i_0} u\|_p } \right) 
\|R_{i_0} u \|_p ,
\end{equation}
where   
$$A(p , \d ) = \frac{p^2 A( \d )}{p-1} \quad {\rm and}\quad 
 A( \d ) \to \infty  \quad {\rm as}\quad   \d \to 0 . $$
For fixed $0 < \a <1  $ a simplified version of  \eqref{19juni14-8}
is the following
\begin{equation}\label{4nov011-3} 
\|W^{(\varepsilon)} (u)\|_p \le \frac{C(p ,  \a,\d )}{1-\a} 
\| u \|_p^{1-\a}
\|R_{i_0} u \|_p^{\a} .
\end{equation}
Specializing \eqref{19juni14-8}, \eqref{4nov011-4} and \eqref{4nov011-3} to the case  $p = 2 $  gives back
  \eqref{19juni14-7}, \eqref{4nov011-1} and \eqref{4nov011-2}.

\paragraph{Review of  \cite{smue2, LMM}.} 
We next  compare the  inequalities 
\eqref{4nov011-4} and \eqref{4nov011-3}
 to the interpolatory estimates 
for  directional  Haar projections 
\cite{LMM}.
Let $$ \{ h_Q^{(\varepsilon)} :  Q\in \cS, \varepsilon \in \cA \}$$
be  the  isotropic Haar system supported on dyadic cubes. 
(See Section \ref{basic} for the definition.) 
The directional Haar projection is
defined  by
$$P^{(\varepsilon)} (u) = \sum_{ Q\in \cS} 
\la u ,  h_Q^{(\varepsilon)}\ra  h_Q^{(\varepsilon)} |Q|^{-1}, \quad\quad u \in  L^2(\bR^n). $$
In \cite{LMM} we proved that  for $ 1 < p < \infty $
and $ \t_p =   \max \{1/2,1/p\} , $ 
\begin{equation}\label{4nov011-5} 
 \|P^{(\varepsilon)} (u)\|_p \le C( p) 
\| u \|_p^{\t_p}
\|R_{i_0} u \|_p^{1-\t_p}, 
\end{equation}
when  $ \varepsilon = (\varepsilon _1, \dots ,\varepsilon_n ) $ and
$\varepsilon_{i_0} =1 .$
Comparing \eqref{4nov011-5}  to 
\eqref{4nov011-3} 
we observe that:
\begin{enumerate}
\item In \eqref{4nov011-3}  the interpolation  exponents for $W^{(\varepsilon)}$  depend
  just on the order $\a$ of the H\"older estimates, and not on the
  value of $p$. 
\item By contrast in \eqref{4nov011-5}, the  $P^{(\varepsilon)}$ estimates show a critical
  transition at $ p = 2 .$ The  exponents in
  \eqref{4nov011-5} are  best
  possible, as shown in  \cite{LMM} Section 8.
Hence  \eqref{4nov011-5} does not 
  arise as the limit of $ \a \to 0 $ from 
   the estimates \eqref{4nov011-3}.
\item As   H\"older estimates are not available for the Haar system
we exploit in \cite{LMM} 
that  the discontinuities of  Haar functions are concentrated at an 
$(n-1)$ dimensional set, and that 
\begin{equation}\label{6nov011-1}  
\int_{\bR^n} | h_Q^{(\varepsilon)} ( x-s ) - h_Q^{(\varepsilon)}(x)
|^p dx \le C |s| \cdot |Q|^{\frac{n-1}{n}}  .
\end{equation} 
\item
It 
remains an open problem to prove interpolatory estimates 
for directional  projections  $W^{(\varepsilon)}$ when the underlying
wavelets satisfy decay estimates {\em only.} 
The particular interest in this question comes from theorems of 
G. Gripenberg \cite{gripen} and P. Wojtaszczyk \cite{woj99} who proved
 that wavelets with decay
\eqref{structure011a} 
form  an unconditional basis in $L^p, ( 1 < p < \infty) $  -- 
{\em  without}  using assumptions on smoothness.  
 
\end{enumerate}
\paragraph{Outlook.} 
In \cite{smue2, LMM}, proving  the weak  semi-continuity  of 
separately convex functionals --as conjectured
by   J. Ball and F. Murat \cite{ballmurat}
 and L. Tartar \cite{tat7} -- provided  the initial motivation for
 estimating Haar projections 
in terms of Riesz transforms. 
For further motivation we refer to the  
analysis of Sverak's counter-examples to
quasi-convexity  in
\cite{smuequasi}.

In the course of development \cite{smue2, LMM}.
the inequalities  \eqref{4nov011-5} gave rise to  
 general questions of  {\em ordering} 
singular integral operators on a given space  by means of  interpolatory estimates.
This includes  the following problems:
\begin{enumerate}
\item Determination of the best possible exponents in interpolatory
  estimates. See \cite{LMM} Section 8 for the sharp exponents between 
Haar projections and Riesz transforms.
\item Extensions to vector valued  singular integral operators.
 R. Lechner \cite{lechner}  obtained the UMD version of \cite{smue2, LMM}

\item Presently  interpolatory estimates between singular integral
  operators 
are known only for the setting of $\bR^n . $ For singular integrals
over 
non commutative groups (e.g.  Heisenberg group, homogeneous  Lie
groups)  such  estimates are open.
See  M. Christ \cite{christ, christ88},
M. Christ and  D. Geller \cite{christgeller}, 
P. G. Lemarie \cite{lemarie}, Folland and Stein \cite{follstei}. 
\end{enumerate}

The results of the present paper
and \cite{smue2, LMM}  are the  first steps in this direction.
 
\section{The Main Results}

Theorem \ref{th1} is the main result
of this paper. The partial coercivity of 
Riesz transforms \eqref{himaug14}  follows immediately from 
Theorem \ref{th1}.

\begin{theor} \label{th1}
Let $1 < p < \infty ,$  
 $ 1 \le i_0 \le n , $ 
and 
$ 
 \varepsilon = (\varepsilon_1, \dots  \varepsilon_n) \in \cA $ with  
$ \varepsilon_{i_0} = 1 .$
If  $ 0 < \a < 1 ,$ then for any $ u \in L^p ( \bR^n ) , $
\begin{equation}\label{19juni141} 
\|W^{(\varepsilon)}(u)\|_p \le 
C \left( 
\|R_{i_0}\|_p^{-\a} \|u\|_p^{1-\a}\|R_ {i_0}(u)\|_p^{\a} + 
 \frac{\|R_{i_0}\|_p ^{1-\a} \|u\|_p^{1-\a} \|R_{i_0} (u)\|_p^\a    - \|R_{i_0} (u)\|_p }{ 2^{1 -\a }   - 1}   \right),
\end{equation}
where   $\|R_{i_0}\|_p = \| R_{i_0} : L^p( \bR^n) \to 
L^p( \bR^n) \| ,$ and  
$$C = C(p ,  \a,\d ) = \frac{p^2 C(  \a,\d )}{p-1} \quad {\rm and}\quad 
 C(  \a,\d ) \to \infty  \quad {\rm as}\quad    \a \to 0   \quad {\rm
   or}\quad    \d \to 0 . $$
If  $ \a = 1 , $ then 
\begin{equation}\label{19juni142} 
 \|W^{(\varepsilon)}(u)\|_p \le C(p , 1,  \d  )  
\left( \|R_{i_0}\|_p^{-1}  \|R_{i_0} ( u) \|_p    + \log\left(\frac{\|u\|_p \|R_{i_0}\|_p}{\|R_{i_0}
    (u)\|_p} \right) \|R_ {i_0}(u)\|_p  \right) .
\end{equation}
The estimate \eqref{19juni142} appears as the limit of  \eqref{19juni141} as  $\a \to 1. $
\end{theor}

\paragraph{Remark.} Clearly \eqref{19juni141} implies that 
$$\|W ^{(\varepsilon)}(u)\|_p \le   \frac{C}{1-\a}
 \|u\|_p^{1-\a}\|R_ {i_0}(u)\|_p^{\a} , $$
and \eqref{19juni142} yields
$$
\|W ^{(\varepsilon)}(u)\|_p \le  C
\log\left( 1 +  \frac{\|u\|_p }{\|R_{i_0}
    (u)\|_p} \right) \|R_ {i_0}(u)\|_p  .$$

\subsubsection*{Partial coercivity of Riesz transforms.}
Theorem \ref{th1} implies  partial-coercivity estimates 
for  
 Riesz transforms. 
  On the    closure of $W^{(\varepsilon)}
 ( L^p ( \bR^n ) )  , $ the Riesz transform $R_{i_0} $ is invertible, provided that $ \varepsilon = 
(\varepsilon _1, \dots ,\varepsilon_n ) $ and
$\varepsilon_{i_0} =1 .$ Indeed, since 
$W^{(\varepsilon)}$   is a  projection,
Theorem \ref{th1} gives 
\begin{equation}\label{himaug14}  
   \|v \|_p \le  C(p ,  \a,\d )\|R_{i_0} v  \|_p ,  \quad v \in W^{(\varepsilon)}
 ( L^p ( \bR^n ) ) .
\end{equation} 
By \eqref{4nov011-5}
the same holds when  $W^{(\varepsilon)}$ is replaced by
$P^{(\varepsilon)}.$
For the concept and background we refer to  T. Kato \cite{kato75},
F. Murat \cite{mur3}, B. Dacorogna \cite{daco1}. 
The interpretation of Theorem \ref{th1} as a 
partial coercivity estimate for Riesz transforms 
\eqref{himaug14} emphasizes  the connection to  \cite{christ88}.

\subsubsection*{The Outline of the Proof.}
We use the pattern of reduction 
applied previously in \cite{smue2, LMM}. In the present paper we 
exploit properties of the discrete Calderon reproducing formula 
going back to Frazier and Jawerth \cite{fraja}.
We start the proof of 
Theorem~\ref{th1} with  a multi-scale analysis of  $ W^{(\varepsilon)} $
using a discrete  Calderon reproducing formula. See \cite{fraja}.  We fix 
$v , w \in C^\infty ( \bR^n) $ so that $ \supp \cF v  , \supp \cF w \sbe [1/2, 2 ] $ and 
$$ 1 = \sum _{\ell \in \bZ }  (\cF v) (2^{\ell}\zeta)   (\cF w) (2^{\ell}\zeta)  .$$
For any multi- index $ \g \in \bN ^n $ and $ N\in \bN $ there exists 
$A = A( \g , N )$ so that
$$ |\pa _\g v(x) | +  |\pa _\g w(x) |  \le A( 1 +|x|) ^{-N}  . $$
Put $v_\ell (x) = 2^{\ell n} v (2^{\ell}x) ,   w_\ell (x) = 2^{\ell n} w(2^{\ell}x) , $  
and form   the convolution product 
 $$ d_\ell(x) = v_\ell*w_\ell(x) , \quad  \ell \in \bZ. $$
For any multi- index $ \g \in \bN ^n $ and $ N\in \bN $ there exists 
$A = A( \g , N )$ so that
\begin{equation}
\label{18juni141}
 |\pa_\g  d_{\ell} (x) |  \le A 2^{(n+|\g|) (\ell)}( 1 +2^{\ell}|x|) ^{-N}  .
\end{equation}
Finally we put  
\begin{equation}
\label{12okt072}
\Delta_{\ell}(u) =  u*d_{\ell} .
\end{equation}
Then as obtained by Frazier Jawerth \cite{fraja}
$$u=\sum^\infty_{\ell=-\infty}\Delta_\ell(u),\quad u \in L^p( \bR^n) ,$$
where convergence 
holds  in $L^p( \bR^n) . $

We denoted by  $\cS $ the collection of all dyadic cubes in $\bR^n . $
Let  $j \in \bZ $ and consider the following subcollection of  $\cS ,$
\begin{equation} \label{25jan068}
\cS_j = \{Q\in \cS :|Q|=2^{-nj}\}.
\end{equation}
The cubes in $\cS_j$ are  pairwise disjoint. 
To   $\ell\in\bZ, $ $\varepsilon \in \cA ,$  and $Q\in \cS_j ,$ define 
\begin{equation}
\label{28mai1}
 f^{(\varepsilon)}_{Q,\ell} = \Delta _{j+\ell} (
 \vp_Q^{(\varepsilon)}), 
\end{equation}
and
\begin{equation}
\label{27mai7}
T_\ell^{(\varepsilon)} (u) = \sum_{Q \in \cS } \la u ,f^{(\varepsilon)}_{Q,\ell}\ra 
 \vp_Q^{(\varepsilon)} |Q|^{-1} .
\end{equation}
Since  the convolution operator  $\Delta_{j+\ell}$ is   self adjoint,  
we arrive at our basic Littlewood-Paley decomposition for the
directional wavelet projection 
$$W^{(\varepsilon)}(u)=\sum^\infty_{\ell=-\infty}T_\ell^{(\varepsilon)} (u).$$

Let  $1 \le i_0 \le n $ 
and 
$ 
\cA_{i_0} = \{ \varepsilon \in \cA : 
 \varepsilon = (\varepsilon_1, \dots  ,\varepsilon_n) 
\quad\text{and}\quad \varepsilon_{i_0} = 1 
\}. $
As observed in \cite{smue2,LMM}, for   $\varepsilon \in \cA_{i_0}$ 
we get 
 $$
T_\ell^{(\varepsilon)}  R^{-1}_{i_0}=T_\ell^{(\varepsilon)}   R_{i_0} + \sum^n_{\substack{i=1\\i \ne i_0}   }
T_\ell^{(\varepsilon)} \bE_{{i_0}}\pa{_i}R_i ,$$
where $R_i $ denotes the $i-$th Riesz transform,
 $\pa{_i}$ denotes the differentiation with respect to the 
$x_i $ variable and $\bE_{{i_0} }$ the 
 integration  with respect to the $x_{i_0}-th  $ 
 coordinate. See \eqref{int121}.  
Hence, putting
\begin{equation}
\label{27mai2}
k_Q^{( \ell, i )} = 
\Delta_{j+\ell}\left( \bE_{{i_0}}\pa{_i}\vp_Q^{( \varepsilon )}\right) ,\quad\quad Q \in \cS_j , 
\end{equation}
we obtain the representation
\begin{equation}
\label{21maerz1}
T_\ell ^{(\varepsilon)} R^{-1}_{i_0} (u) = T_\ell^{(\varepsilon)}
  R_{i_0} (u) +
\sum_{Q \in \cS}  \sum^n_{\substack{i=1\\i \ne i_0}   }
\la R_i (u) , k_Q^{( \ell, i )}\ra \vp_Q^{( \varepsilon )} |Q|^{-1}.
\end{equation} 

The following two theorems record the norm estimates for the 
operators  $T_{\ell}^{(\varepsilon)}$ and $T_\ell^{(\varepsilon)}
R^{-1}_{i_0}$ 
by which we obtain Theorem~\ref{th1}. 
First we treat the case $ \ell > 0.$ It  displays the 
crucial dependence  on  the H\"older
  exponent of the admissible wavelet system. 
Below and throughout the paper the constants $C(p ,  \a,\d ) >0$
satisfy the conditions
 $$C(p ,  \a,\d ) =  \frac{p^2 C(  \a,\d) }{p-1}, \quad\text{where}\quad 
 C(  \a,\d ) \to \infty, \quad\text{ as }\quad \a \to 0 , \quad\text{ or } \quad\d \to 0 .$$

\begin{theor}\label{th2a} Let $ \d > 0$ and  $ 0 < \a \le 1 $ be the 
decay and H\"older
  exponents of the admissible wavelet system specified in  
\eqref{structure011}. 

Let $1 < p  < \infty ,$ 
 $\ell \ge 0 $ 
and  $\varepsilon \in\cA .$ Then  $T^{(\varepsilon)}_{\ell}$ satisfies the
norm estimates, 
\begin{equation}\label{8jan3}
\| T_{\ell}^{(\varepsilon)} \|_p \le 
C(p,\a,\d) 2^{-\ell \a}  .
\end{equation}

Let $1 \le i_0 \le n, $ and $\varepsilon \in \cA_{i_0}$ then  
\begin{equation}\label{81jan3}
\| T_{\ell} ^{(\varepsilon)}R_{i_0}^{-1} \|_p   \le 
C(p,\a,\d) 2^{\ell -\ell\a}  .
\end{equation}
\end{theor}

When $ \ell < 0 $ we get  exponents
  independent of  the H\"older condition.
\begin{theor}\label{th2b} Let $ \d > 0$ and  $ 0 < \a \le 1 $ be the 
decay and H\"older
  exponents of the admissible wavelet system specified in  
\eqref{structure011}. 
Let $1 < p < \infty .$  
Let $\ell \le 0 .$ Then for $\varepsilon \in \cA $
the operator $T_{\ell} ^{(\varepsilon)}$ satisfies the
norm estimates, 
\begin{equation}\label{82jan3}
\| T_{\ell} ^{(\varepsilon)} \|_p   \le
C (p, \a, \d) 2^{-|\ell |} |\ell|. 
\end{equation}
If  moreover $1 \le i_0 \le n, $ and $ \varepsilon \in \cA_{i_0} ,$
then  
\begin{equation}\label{82jan3xx}
\| T_{\ell} ^{(\varepsilon)} R_{i_0}^{-1} \|_p \le 
C(p,\a,\d) 2^{-|\ell |} |\ell | . 
\end{equation}

\end{theor}
\subsubsection*{The Proof of Theorem \ref{th1}.}
 Theorem~\ref{th2a} and  Theorem~\ref{th2b} yield the 
proof of Theorem \ref{th1} as follows.
Fix $ u \in L^p( \bR^n ) .$ Define $M \in \bN$  by  the relation 
\begin{equation}
\label{17feb3}
2^{M-1} \le \frac{\|u\|_p \|R_{i_0}\|_p}{\|R_{i_0} (u)\|_p}\le 2^{M} .
\end{equation}
First we fix the H\"older exponent wavelet system as  $ 0 < \a < 1 . $  
Thereafter we consider the limit as
$ \a \to 1 . $ By 
Theorem~\ref{th2a} and  Theorem~\ref{th2b}
 there exists $C =  C(p,\a,\d) $ so that 
$$\sum^\infty_{\ell=M}\| T_{\ell} ^{(\varepsilon)} \|_p \le C  2^{-M\a}, $$
and
$$\sum^{M-1}_{\ell= -\infty}\| T_{\ell} ^{(\varepsilon)} R_{i_0}^{-1}
\|_p \le C  \left(    \sum^{M-1}_{\ell= 0}  2^{\ell -\a \ell} \right)
= C       \frac{ 2^{M -\a M}   - 1 }{ 2^{1 -\a }   - 1}.  $$
The constants  $C(p,\a,\d)$ stay bounded as $  \a \to 1 . $ 
Since 
$
W ^{(\varepsilon)}(u)=\sum^\infty_{\ell=-\infty}T_\ell ^{(\varepsilon)} (u) 
$
 triangle inequality gives 
\begin{equation}\label{zerl11}
\begin{aligned}
 \|W ^{(\varepsilon)}(u)\|_p &\le  \sum ^{\infty}_{\ell=M}\|T_\ell ^{(\varepsilon)}\|_p\|
  u\|_p+\sum^{M-1}_{\ell= -\infty}\|T_\ell ^{(\varepsilon)} R^{-1}_{i_0}\|_p\, \|R_{i_0}(u)\|_p\\
&\le   C \left(   2^{-M\a}\|u\|_p +  
  \frac{ 2^{M -\a M}   - 1 }{ 2^{1 -\a }   - 1} \|R_{i_0}(u)\|_p \right) .
\end{aligned}
\end{equation}
 Inserting the value of $M$ specified in \eqref{17feb3} gives
the following upper bound for  \eqref{zerl11},
\begin{equation}\label{18juni142}
C \left( 
\|R_{i_0}\|_p^{-\a} \|u\|_p^{1-\a}\|R_ {i_0}(u)\|_p^{\a} + 
 \frac{\|R_{i_0}\|_p ^{1-\a} \|u\|_p^{1-\a} \|R_{i_0} (u)\|_p^\a    - \|R_{i_0} (u)\|_p }{ 2^{1 -\a }   - 1}   \right).
\end{equation}
The term arising in \eqref{18juni142}  has a well defined limit as  
   $\a \to 1. $ Indeed, by  l'H\^{o}pital's rule,
$$
\lim_{\a \to 1 }  \frac{ 2^{M -\a M}   - 1 }{ 2^{1 -\a }   - 1} = M  , $$
and hence 
\begin{equation}\label{18juni143}
\lim _{\a \to 1} \eqref{18juni142} =  {C (p,\d) } \left( \|R_{i_0}\|_p^{-1}  \|R_{i_0} ( u) \|_p    + \log\left( \frac{\|u\|_p \|R_{i_0}\|_p}{\|R_{i_0}
    (u)\|_p} \right) \|R_ {i_0}(u)\|_p  \right) . \end{equation}
\endproof
\paragraph{Remark 1.} The complicated form of \eqref{18juni142} was used to obtain the 
limit estimate \eqref{18juni143}. 
For fixed $ \a < 1 $  we  we may simplify  the upper bound \eqref{18juni142}
as follows  
\begin{equation}\label{18juni144}
\|W ^{(\varepsilon)}(u)\|_p \le   \frac{C}{{1-\a}}
 \|u\|_p^{1-\a}\|R_ {i_0}(u)\|_p^{\a} ;
\quad\quad
 \quad 1 < p < \infty . 
\end{equation} 
Also  \eqref{18juni143} may be simplified further,
\begin{equation}\label{18juni1444}
\|W ^{(\varepsilon)}(u)\|_p \le  C
\log\left( 1 +  \frac{\|u\|_p }{\|R_{i_0}
    (u)\|_p} \right) \|R_ {i_0}(u)\|_p  .
\quad\quad
 \quad 1 < p < \infty . 
\end{equation} 
\paragraph{Remark 2.}  An alterative  proof   of  \eqref{18juni143} 
may be deduced  directly from 
Theorem~\ref{th2a} and  Theorem~\ref{th2b}. 
 Define $M \in \bN$  by   
\eqref{17feb3}. Then 
$$\sum^\infty_{\ell=M}\| T_{\ell} ^{(\varepsilon)} \|_p \le C(p, \d)  2^{-M}
\quad\text{and}\quad
\sum^{M-1}_{\ell= -\infty}\| T_{\ell} ^{(\varepsilon)} R_{i_0}^{-1}
\|_p \le C(p, \d) M^.
$$ 
Hence, by the  triangle inequality,  
\begin{equation}\label{zerl111}
\begin{aligned}
 \|W ^{(\varepsilon)}(u)\|_p \le   C(p, \d)   2^{-M} \|u\|_p +   C(p,\d) M \|R_{i_0}(u)\|_p.
\end{aligned}
\end{equation}
 As $M$ is given  by \eqref{17feb3} 
 the right hand side of \eqref{zerl111} is dominated by 
$$ {C (p,\d) } \left( \|R_{i_0}\|_p^{-1}  \|R_{i_0} ( u) \|_p    + \log\left( \frac{\|u\|_p \|R_{i_0}\|_p}{\|R_{i_0}
    (u)\|_p} \right) \|R_ {i_0}(u)\|_p  \right);
\quad\quad
 \quad 1 < p < \infty . $$
\endproof

\paragraph{The organization of the paper:} In  Section
\ref{technical} we prove point-wise estimates for the decay and
smoothness of the systems $\{f^{(\varepsilon)}_{Q,\ell}\} $
and  $\{k_Q^{( \ell, i )}
\} $ defined in \eqref{28mai1} and \eqref{27mai2}. 
The cases $ \ell > 0 $ and $\ell \le 0 $ are given different
treatment.

In Section \ref{basic}
we present two general tools used to reduce estimates for 
 integral operators to those of rearrangements.

In Section \ref{proof} we combine the preparatory theorems of 
 Section \ref{technical} and Section \ref{basic}
to  prove Theorem~\ref{th2a}.
Section \ref{simple}  contains the proof of
Theorem \ref{th2b}.
\section{Wavelets and Convolution}

\label{technical}

Our basic concern are the norm estimates for the operators
$ T_{\ell}^{(\varepsilon)}$ and  $  T_{\ell} ^{(\varepsilon)}
R_{i_0}^{-1}$
as formulated in 
 Theorem~\ref{th2a} and  Theorem~\ref{th2b}.
We showed in the introduction that this amounts to proving estimates for 
operators 
$$
X (u) = \sum_{Q \in \cS}  \la u , k_Q^{( \ell, i )}\ra \vp_Q^{( \varepsilon )} |Q|^{-1},
\quad\text{ where}\quad 
k_Q^{( \ell, i )} = 
\Delta_{j+\ell}( \bE_{{i_0}}\pa{_i}\vp_Q^{( \varepsilon
    )}), 
\quad Q \in \cS_j , \, \varepsilon \in \cA_{i_0}, \, i \ne i_0 , 
$$
and 
$$
Y (u) = \sum_{Q \in \cS}  \la u ,  f^{(\varepsilon)}_{Q,\ell} \ra \vp_Q^{( \varepsilon )} |Q|^{-1},
\quad\quad\text{where}\quad  
 f^{(\varepsilon)}_{Q,\ell} = \Delta _{j+\ell} (
 \vp_Q^{(\varepsilon)}),\quad Q \in \cS_j  . 
$$
Recall that in  \eqref{12okt072} we defined  the operator $ \Delta_{j+\ell} $ 
as convolution with  $d_{j+\ell} ,$
where for any multi- index $ \g \in \bN ^n $ and $ N\in \bN $ there exists 
$A = A( \g , N )$ so that
\begin{equation}
\label{12okt075}
 |\pa_\g  d_{j+\ell} (x) |  \le A 2^{(n+|\g|) (j+\ell)}( 1 +2^{j+\ell}|x|) ^{-N}  .
\end{equation}
In Lemma \ref{basiclemma} through Lemma \ref{9dez11a}
 we record  the point-wise estimates for the 
systems  $\{f^{(\varepsilon)}_{Q,\ell} \}$ and $\{k_Q^{( \ell, i )}\}
 $  as  needed for the purpose of this paper.
Those are the basis for the norm inequalities of the operators $X$
and $Y $ defined above.

\subsection{Point-wise estimates for  $\Delta _{j+\ell} (
  \vp_Q^{(\varepsilon)}). $ }

The following Lemma records basic point-wise estimates for
 $f^{(\varepsilon)}_{Q,\ell},$ $\ell \ge 0$ and its gradient.

\begin{lemma}\label{basiclemma}
Assume that $\{ \vp_Q^{(\varepsilon)} :\,Q \in \cS \} $ satisfies
\eqref{structure011}. The system
 $\{f^{(\varepsilon)}_{Q,\ell}:\,Q \in \cS ,\ell \ge 0\} $  defined by
 \eqref{28mai1} 
 satisfies  these basic estimates:
\begin{subequations} \label{27mai6}
\begin{equation}
\label{27mai6z}
| f^{(\varepsilon)}_{Q,\ell}(x)| \le C
2^{-\a\ell} \left( 1 + \frac{\dist(x, Q) }{s(Q) }
 \right)^{-n(1+\d)},  \quad\quad x \in \bR^n .\end{equation}  
\begin{equation}\label{27mai6c}
|\nabla  f^{(\varepsilon)}_{Q,\ell} (x)| \le  C 2^{-\a\ell} 2^{j+\ell} \left( 1 + \frac{\dist(x, Q) }{s(Q) }
\right)^{-n(1+\d)} , \quad\quad x \in \bR^n .
\end{equation}
\begin{equation}\label{27mai6zero}
\int_{\bR^n}f^{(\varepsilon)}_{Q,\ell} (x)dx = 0 .
\end{equation}
\end{subequations}
\end{lemma}
\proof
Let $ x \in \bR^n ,  $ $ Q \in \cS_j $ and $ \ell > 0 . $
Let ${A_x} =  \{ t : | x- t | \le C
2^{-j -\ell } \} . $
Since $\int d_{j + \ell} (x -t) dt = 0 $ and $ t \to  d_{j + \ell} (x -t) $ is 
centered at $A_x$
we get: 
$$ 
\begin{aligned} 
| d_{j + \ell} *  \vp_Q^{(\varepsilon)} (x)| & = \left |\int_{\bR^n }   d_{j +
  \ell} ( x-t)( \vp_Q^{(\varepsilon)} (t) - \vp_Q^{(\varepsilon)} (
x)) dt \right|\\
& \le \int_{\bR^n }
 |  d_{j +
  \ell} ( x-t)| \cdot | \vp_Q^{(\varepsilon)} (t) - \vp_Q^{(\varepsilon)} (
x)| dt 
\\
& \le C\int_{\bR^n }
| d_{j +
  \ell} ( x-t)| dt   \diam (A_x) ^{\a}  s (Q)^{-\a} 
\left( 1 + \frac{\dist(x, Q) }{s(Q) }
 \right)^{-n(1+\d)} 
. 
\end{aligned}
$$
Invoking that $ \diam (A_x) \le C 2^{-j-\ell },$   $s (Q) = 2^{-j }$
and  $\int_{\bR^n}
 |  d_{j +
  \ell} ( x-t)|dt \le C $ yields \eqref{27mai6z}.

In a similar fashion we obtain the remaining estimates
\eqref{27mai6c}.
Put $ \widetilde   d_{j +\ell} = 2 ^{-(j+\ell) } \nabla d_{j+\ell} . $
Repeating the above argument  with  $ d_{j +\ell}$ replaced by 
$\widetilde   d_{j +\ell} $ we get
$$ 
\begin{aligned} 
|\nabla d_{j + \ell} *  \vp_Q^{(\varepsilon)} (x)| & = 2^{j+\ell}
\left |\int_{\bR^n }   \widetilde d_{j +
  \ell} ( x-t)( \vp_Q^{(\varepsilon)} (t) - \vp_Q^{(\varepsilon)} (
x)) dt \right|\\
& \le 2^{j+\ell} \int_{\bR^n }
 |  \widetilde d_{j +
  \ell} ( x-t)| \cdot | \vp_Q^{(\varepsilon)} (t) - \vp_Q^{(\varepsilon)} (
x)| dt \\
& \le C2^{j+\ell} \int_{\bR^n }
| \widetilde  d_{j +
  \ell} ( x-t)|  dt  \diam (A_x) ^{\a}  s (Q)^{-\a} 
\left( 1 + \frac{\dist(x, Q) }{s(Q) }
 \right)^{-n(1+\d)} \\
&\le C 2^{-\a \ell} 2^{j+\ell}\left( 1 + \frac{\dist(x, Q) }{s(Q) }
 \right)^{-n(1+\d)} . 
\end{aligned}
$$

\endproof

\subsubsection*{Compactly supported Wavelets.}
Recall compactly supported wavelets.
Let  
 $\{ \psi_K^{(\b)} /\sqrt{| K|}
 : K \in \cS, \b \in \cA 
\} $
be an orthonormal basis in $L^2(\bR^n),$
satisfying  $\int \psi_K^{(\b)} = 0 $ 
and
the following structure conditions,
\begin{equation}\label{13april_1}
 \supp\psi_K^{(\b)} \sbe C\cdot K,
\quad \quad |\psi_K^{(\b)}| \le C ,
\quad\quad \Lip(\psi_K^{(\b)}) \le C s(K)^{-1}.
\end{equation}
We often  write 
in place  of  $\{\psi_K^{(\b)}\}$ 
just $\{ \psi_K \}.$  The existence of compactly supported wavelets
was proven by I. Daubechies, see \cite{coifme}.
\subsubsection*{Low frequency slices of  $\Delta _{j+\ell} (
  \vp_Q^{(\varepsilon)}). $ }
Here we prove point-wise estimates for decay and regularity of the 
low frequency slices of 
$
 f^{(\varepsilon)}_{Q,\ell} 
$ 
when  $\ell \ge 0. $
We define those slices using a compactly supported wavelet basis $\{
\psi_K \}$ satisfying \eqref{13april_1}.
Fix 
 $ k \in \bZ \sm \bN  $ and define
\begin{equation}\label{12april_1}
p_{Q}=
\sum_{K\in \cS_{j +k}}
 \left \la  f^{(\varepsilon)}_{Q,\ell},\psi_K
\right\ra \psi_K |K|^{-1} , \quad   Q \in\cS_j ,  \quad j \in \bZ .
\end{equation}
Note that  for  $ k \in \bZ \sm \bN . $
and  $ Q \in\cS_j ,  j \in \bZ $  
there exists a unique cube $K_0 = K_0(Q) $ so that 
\begin{equation}\label{12april_2}
K_0 \spe Q \, , \quad K_0 \in \cS_{j+k} .
\end{equation}
Pointwise estimates for  $ p_{Q}$ and    $\nabla p_{Q} $ are as follows.
\begin{lemma} \label{12april_3} Let  $ \ell \in \bN $ and $ k \in \bZ \sm \bN . $
Let $ \g_0 = \min\{ \d n/2 , 1 \} .$ 
Then the  system of slices 
defined by \eqref{12april_1} satisfies the following estimates:
\begin{subequations}\label{12april_4A} 
\begin{equation}\label{12april_4B}
| p_{Q} (t)|  \le C  2^{-\a\ell} 2^{k  (n + \g_o) }
\left(1 +\frac{\dist(t , K_0)}{s(K_0)} \right)^{-n(1+\d)}, 
\end{equation}
\begin{equation}\label{12april_4C}
|\nabla p_{Q} (t)| \le Cs(K_0)^{-1} 2^{-\a\ell}
2^{k  (n + \g_o) }
\left(1 +\frac{\dist(t , K_0)}{s(K_0)} \right)^{-n(1+\d)}. 
\end{equation}
\end{subequations}
where $Q \in \cS $ and $K_0 = K_0(Q) $ is defined by 
\eqref{12april_2}.
\end{lemma}
\proof
Fix 
 a dyadic cube $Q \in \cS_j , $  $j \in \bZ .$
Determine $K_0 \in \cS_{j +k} $  so that $ Q \sbe K_0 .$
For any $ \mu  \in \bZ^n $ and   $K = K_0 + \mu \cdot s(K_0 ) $ we prove the
following coefficient estimate
\begin{equation}\label{12april_5}
  |\la  f^{(\varepsilon)}_{Q,\ell}  ,\psi_K \ra | |K|^{-1} 
\le  2^{-\a\ell}  2^{-|k| (n+\g_0) } (1+|\mu |)^{-n(1+\d)} . 
\end{equation}
Consider first the case $ |\mu | \ge 4 . $ Then 
Lemma \ref{basiclemma} gives
\begin{equation}\label{12april_6}
  |\la  f^{(\varepsilon)}_{Q,\ell},\psi_K
\ra | |K|^{-1} \le  2^{-\a\ell} |\mu |^{-n(1+\d)} 2^{k n(1+\d) }. 
\end{equation}
Note that   \eqref{12april_6} implies \eqref{12april_5} by arithmetic.

Next consider   $ |\mu | \le 4 .$ 
We use that $ f^{(\varepsilon)}_{Q,\ell}$ is of vanishing
mean and rewrite  
\begin{equation}\label{12april_7}
\la  f^{(\varepsilon)}_{Q,\ell}  ,\psi_K \ra =
\int_{\bR^n}  f^{(\varepsilon)}_{Q,\ell} (t)(\psi_K (t) -
\psi_K (t_Q) ) dt  \end{equation}
where $t_Q \in Q . $ 
We decompose the domain of integration as follows.
Let  $A_0(Q) = Q $ and $ A_i ( Q ) = 2^i\cdot Q \sm 2^{i-1}\cdot Q $
where $2^i\cdot Q$ is the  cube with side-length $ 2^i s( Q ) $ and the same
center as $Q.$  Thus defined the sets   $ A_i ( Q ) = 2^i\cdot Q \sm
2^{i-1}\cdot Q $
form a decomposition of $ \bR^n .$
Hence 
the right hand side of \eqref{12april_7} 
is bounded by 
\begin{equation}\label{him1}
\sum_{i = 0 }^{|k| } \int_{ A_i ( Q )} 
|  f^{(\varepsilon)}_{Q,\ell} (t)
(\psi_K (t) -
\psi_K (t_Q) )| dt   + \sum_{i = |k|+1  }^{\infty } \int_{ A_i ( Q )} 
|  f^{(\varepsilon)}_{Q,\ell} (t)| dt \end{equation} .
For $ i \le |k| $ we exploit the Lipschitz estimates for $\psi_K$ and use that 
$\diam (A_i (Q) ) \le C 2^i s (Q ) \le C s(K) .$ Thus we get  
$$
| \psi_K (t) -
\psi_K (t_Q) | \le C ( 2^i s (Q) ) ^{\g_0}(\Lip\psi_K) ^{\g_0} ,  \quad\quad t \in  A_i ( Q ).$$
Invoking also the basic estimates of Lemma \ref{basiclemma}
gives 
\begin{equation}\label{9oct121}   \int_{ A_i ( Q )} 
|  f^{(\varepsilon)}_{Q,\ell} (t)
(\psi_K (t) -
\psi_K (t_Q) )| dt \le   2^{-\a\ell}  2^{-i n(1+\d) } ( 2^i s(Q) ) ^{\g_0}(\Lip\psi_K) ^{\g_0} 
| 2^i \cdot Q | .  \end{equation}
We take  the sum over $ i \le |k|   $ in \eqref{9oct121} and  get the following 
upper bound for
the  first sum in \eqref{him1}, 
$$ 2^{-\a\ell}\frac{ (s( Q))^{\g_0} |Q| } {( s( K) )^{\g_0} } \sum_{i = 0 }^{k }
2^{-in\d +i{\g_0}   } 
\le C 2^{-\a\ell}\frac{ (s( Q) )^{n+\g_0}  } {(s( K)  )^{\g_0} } .
$$
On the other hand if $ i \ge |k| $ we have again by  Lemma \ref{basiclemma} that 
$$ \int_{ A_i ( Q )} 
|  f^{(\varepsilon)}_{Q,\ell} (t)
\le C  |Q| 2^{-\a\ell}  2^{-i n\d } .   $$
Since $ |K | 2^{-n|k|} = | Q| ,$  summing over $ i \ge |k| $ gives 
$$
 |Q| \sum_{i = |k|+1  }^{\infty }   2^{-\a\ell}  2^{-i n\d }\le 
  C |K |  2^{-\a\ell}  2^{- |k| n(\d +1)  } . 
$$
Recall now that we are treating the case $|\mu| \le 4 .$ Hence we may
rephrase the above as
$$   |\la  f^{(\varepsilon)}_{Q,\ell}  ,\psi_K \ra | |K|^{-1} 
\le C  2^{-\a\ell}  2^{-|k| (n+\g_0) } (1+|\mu |)^{-n(1+\d)} . $$ 
By the definition of the slices \eqref{12april_1} 
we obtain the point-wise estimates \eqref{12april_4A} from  \eqref{12april_5}.
\endproof
\subsubsection*{The Scalar Products $\la \Delta _{j+\ell} 
( \vp_Q^{( \varepsilon )}) , \Psi_K \ra ,\,\, K \in \cS_{j + k} , 
\,\, 0 \le k \le \ell . $}
Here we record a short but crucial consequene of  Lemma \ref{basiclemma}. It is here where we 
explicitley exploit that our multi-scale analysis   $\{ d_\ell\} $ is based on  
Calderon's  reproducing formula  and admits a   factorization as 
$$ d_\ell= v_\ell * w_\ell . $$

\begin{lemma} \label{scalar-prod} Let $ j \in \bZ , \ell \in \bN  , $ and $ 0 \le k \le \ell . $ For $Q \in \cS_j ,$ 
$ K \in \cS _{j + k } , $
\begin{equation}\label{17juni14} 
| \la   \Delta_{j+\ell} ( \vp_Q^{( \varepsilon )}) , \Psi_K \ra | \le 
C 2^{-\a \ell } 2^{k-\ell} |K | \left( 1 +\frac{\dist (  K , Q )}{ s(Q)} \right) ^{-n(1 + \d)} . 
\end{equation}
\end{lemma}
\proof
Recall that  $ \Delta_{j+\ell} ( \vp_Q^{( \varepsilon )})  = d_{j+\ell} * \vp_Q^{( \varepsilon )} $
where $  d_{j+\ell} = v_{j+\ell}* w_{j+\ell} . $ Hence 
$$  \la   \Delta_{j+\ell} ( \vp_Q^{( \varepsilon )}) , \Psi_K \ra =
 \la w _{j+\ell} *\vp_Q^{( \varepsilon )} , v_{j+\ell}*\Psi_K \ra . $$
By Lemma \ref{basiclemma} we get
\begin{equation}\label{him2}   |w _{j+\ell} *\vp_Q^{( \varepsilon )}(x ) | \le  2^{-\a \ell } \left( 1 +\frac{\dist (  x , Q )}{ s(Q)} \right) ^{-n(1 + \d)} .\end{equation}
Similarly, using that $ j +\ell \ge j + k \ge j ,$ 
 the proof of  Lemma \ref{basiclemma} gives routinely the estimate, 
\begin{equation}\label{him3}  |v_{j+\ell} *\psi_K (x ) | \le  2^{ k -\ell } \left( 1 +\frac{\dist (  x , K )}{ s(K)} \right) ^{-4n} . \end{equation}
Taking into account that $ s(K ) \le s(Q) $  we get
$$ \int_{\bR^n }  
\left( 1 +\frac{\dist (  x , Q )}{ s(Q)} \right) ^{-n(1 + \d)} \left( 1 +
\frac{\dist (  x , K )}{ s(K)} \right) ^{-4n} dx \le C |K|\left( 1 +\frac{\dist (  K  , Q )}{ s(Q)} \right) ^{-n(1 + \d)} . $$
Combining this with the pointwise estimates \eqref{him2} and \eqref{him3}  gives \eqref{17juni14}. 
\endproof
  
\subsection{Point-wise estimates for $\Delta_{j+\ell}\left( \bE_{{i_0}}\pa{_i}\vp_Q^{( \varepsilon )}\right).$  }

We turn to  the analysis of the system
 $k_Q^{( \ell, i )} = \Delta_{j+\ell}\left( \bE_{{i_0}}\pa{_i}\vp_Q^{( \varepsilon )}\right) 
 $ as defined by \eqref{27mai2}. The
 cases  $\ell \ge 0 $ and  $\ell \le 0 $ will be treated separately. 
We begin with the case  $\ell \ge 0 .$

\begin{lemma} \label{basiclemma2} Let $\varepsilon \in \cA_{i_0} .$
The system  $\{k_Q^{( \ell, i )}:\,Q \in \cS ,\, i \ne i_0,\, \ell \ge 0\} $
defined by \eqref{27mai2} satisfies the structural conditions,
\begin{subequations}\label{27mai3} 
\begin{equation}
\label{27mai3a} 
| k_Q^{( \ell, i )} (x) | \le C 2^{\ell-\a \ell } 
\left( 1 +\frac{\dist (  x , Q )}{ s(Q)} \right) ^{-n(1 + \d)}  , \end{equation}
\begin{equation}\label{27mai3b}
|\nabla  k_Q^{( \ell, i )} (x) | \le C 2^{j+\ell}  2^{\ell -\a \ell } 
\left( 1 +\frac{\dist (  x , Q )}{ s(Q)} \right) ^{-n(1 + \d)} 
\end{equation}
\begin{equation}\label{27mai3c} 
\int_{\bR^n}  k_Q^{( \ell, i )} (x) dx = 0.\end{equation}
\end{subequations}
with $C>0$ independent of $ Q\in \cS ,\,$  $i \ne i_0,\,$ or $ \ell \ge 0 .\, $
\end{lemma}
\proof 
Fix $ Q \in \cS , $ and $ x \in \bR^n . $
Put $ e_Q = \bE_{{i_0}}\vp_Q^{( \varepsilon )} , $
then 
$$\begin{aligned} 
k_Q^{( \ell, i )} (x) 
&= \int_{\bR^n}    \pa _i d_{j + \ell} (x-t)e_Q ( t ) dt .
\end{aligned} 
$$ 
By  \eqref{structure011} for admissible wavelets we get  
 \begin{equation}\label{4nov011-a}
| e_Q  (x) | \le C s(Q) 
\left( 1 +\frac{\dist (  x , Q )}{ s(Q)} \right) ^{-n(1 +
  \d)}, \quad\quad x \in \bR^n ,\end{equation}
and
 \begin{equation}\label{4nov011b}
| e_Q  (t) -  e_Q  (s) | \le C |s-t|^\a s(Q)^{1-\a} 
\left( 1 +\frac{\dist (  x , Q )}{ s(Q)} \right) ^{-n(1 + \d)}, 
\quad\quad |s-t | \le s(Q) .
\end{equation}
Hence, with  \eqref{12okt075} and \eqref{4nov011b}, we obtain
\eqref{27mai3a} 
and \eqref{27mai3b} by 
repeating the proof
of Lemma \ref{basiclemma}.
It remains to check  \eqref{27mai3c}.
Since  $ \Delta_{j+\ell} $  commutes with 
differentiation, 
$k_Q^{( \ell, i )} = 
\pa{_i}\Delta_{j+\ell}(e_Q    ) . $
Hence the decay of  $e_Q$  and $\Delta_{j+\ell} e_Q$  imply 
$$ \int_{\bR^n}k_Q^{( \ell, i )} (x) dx = 0 ,$$
that is   \eqref{27mai3c}.

\endproof
Next  we treat the case $ \ell \le 0 . $
\begin{lemma}\label{9dez11a}
The family 
$\{k_Q^{( \ell, i )}:\,Q \in \cS ,\, i \ne i_0,\, \ell \le 0\} ,$
satisfies these conditions 
\begin{subequations}\label{27mai4}
\begin{equation}
\label{27mai4a}
| k_Q^{( \ell, i )} (x) | \le C  2^{- (n + 1 )  |\ell | } \left( 1 +
2^{- |\ell|}\frac{\dist (  x , Q )}{ s(Q)} \right) ^{-n(1 + \d)}  , 
\end{equation}
\begin{equation}
\label{27mai4b}
|\nabla  k_Q^{( \ell, i )} (x) | \le C 2^{j - |\ell|}  2^{- (n + 1 )
  |\ell | } 
\left( 1 +
2^{- |\ell|} \frac{\dist (  x , Q )}{ s(Q)} \right) ^{-n(1 + \d)}  ,
\end{equation}
\begin{equation}
\label{27mai4c}
\int_{\bR^n}  k_Q^{( \ell, i )} (x) dx = 0.
\end{equation}
\end{subequations}
where  $C>0$ is independent of $ Q\in \cS ,\,$  $ i \ne i_0,$ or $ \ell 
\le 0 .\, $
\end{lemma}
\proof Fix $ Q \in \cS , $ and $ x \in \bR^n . $
Put again $ e_Q =  \bE_{{i_0}}\vp_Q^{( \varepsilon )} $ so that   
%
$$ k_Q^{( \ell, i )} (x)=\int_{\bR^n}  e_Q (x - t ) \pa _i d_{j + \ell} (t) 
dt . $$
By \eqref{structure011} 
we have 
\begin{equation}\label{28ap011-D}
 |e_Q ( x  )| \le C \left( 1 + \dist (  x , Q ) / s(Q) \right) ^{-n(1
  + \d)} s(Q) 
\end{equation}
 By \eqref{12okt075} we get  
\begin{equation}\label{28ap011-E} \|  \pa _i d_{j +\ell} \| _1 \le C 2^{j - |\ell|} , 
\quad\text{and}\quad  
 \|  \pa _i d_{j + \ell} \| _{\infty} \le C  2^{n(j - |\ell|)} 2^{j
  - |\ell|} .\end{equation}
We distinguish  between the following cases:
\begin{enumerate}
\item  $\dist ( x ,Q ) \ge 2^{|\ell | } s(Q) $
\item  $\dist ( x, Q ) \le 2^{|\ell | } s(Q) $
\end{enumerate}

In the first case select $ \mu \in \bZ ^n $ so that 
$ x \in 2^{|\ell | }\cdot Q +  \mu 2^{|\ell |} s (Q) . $ 
Then 
we have with \eqref{28ap011-D} and \eqref{28ap011-E}
\begin{equation}
\begin{aligned}
\int_{\bR^n}  e_Q (x - t ) \pa _i d_{j +\ell} (t) dt 
 & \le  2^{- |\ell| }  |\mu | ^{-n(1+ \d)}  2^{ -n |\ell| ( 1 + \d )} \\
&\le \left( 1 +   2^{- |\ell| } \dist (  x , Q ) / s(Q) \right) ^{-n(1+ \d)} 2^{ -n |\ell| ( 2 + \d )} . 
\end{aligned}
\end{equation}

In the second case we have  $\dist ( x, Q ) \le 2^{|\ell | } s(Q).$
Select $ k_0 $ so that
$$ 2^{k_0} s(Q) \le \dist ( x, Q ) \le 2^{k_0+1} s(Q)  
.$$
We may assume that $ k_0 \ge 1 . $
Define the disk
$$A_0 = \{ y \in \bR^n : | y - x | \le \dist( x ,Q )  \}$$
and the annuli 
$$
A_k  = \{ y \in \bR^n :  2^{k-1} (\dist( x ,Q ) ) \le | y - x |
\le 2^k(\dist( x ,Q ) ) 
\} .$$ 
Use \eqref{28ap011-D} and \eqref{28ap011-E} to obtain
\begin{equation} \label{27ap011-14}
\begin{aligned}
\int_{\bR^n}  |e_Q (x - t ) \pa _i d_{j + \ell} (t)| dt & \le 
\sum_{k = 0 } ^\infty \int_{ A_k \sm  A_{k - 1 } }   |e_Q (x - t ) \pa
_i d_{j +\ell} (t)| dt \\
& \le C\sum_{k = 0 } ^\infty |A_k | 
\left( 1 +   2^{k } \dist (  x , Q ) / s(Q) \right) ^{-n(1
  + \d)} s(Q) \|  \pa _i d_{j +\ell} \| _{\infty}  \\
& \le C \sum_{k = 0 } ^\infty  (2^{k } 2^{k_0 })  ^{-n \d}   2^{-( n +
  1 ) |\ell| } .
\end{aligned}
\end{equation}
Clearly \eqref{27ap011-14} gives
\eqref{27mai4a}. 
The gradient estimates \eqref{27mai4b} follow from  \eqref{12okt075} and  \eqref{27mai4a}.
\endproof

\section{Review of Basic Dyadic Operations}
\label{basic}
In this section we prove two auxiliary results on rearrangement
operators. 
The norm estimates for the operators 
$T_\ell^{(\varepsilon)} , \,  T_\ell^{(\varepsilon)} R_{i_0}^{-1}$
will be obtained as applications of Proposition \ref{19april011} and 
Theorem \ref{11julitheorem1}. 
\subsection{The Haar system }
We recall the definition of the isotropic Haar system and its
equivalence to admissible wavelet systems.
We use  \cite{figsingular}, \cite{coifme}, \cite{MR2157745} as sources.
Let $I$ be a dyadic interval and $h_I $ be the $L^\infty$ normalized 
Haar function supported
on $I.$ Thus $h_I = 1 $ on the left half of $I$ and  $h_I = -1 $
on the right half of $I .$
Given a dyadic cube $ Q = I_1 \times \cdots \times I_n  $
and a direction $\varepsilon = (\varepsilon_1 , \dots , \varepsilon _n
) \in \cA $
we define the Haar function 
$$
 h_Q^{(\varepsilon)} ( x ) = \prod_{i = 1 }^n
 h_{I_i}^{\varepsilon_i}(x_i) ,
\quad  x = (x_1, \dots , x_n ) .$$ 
 The Haar system
$ \{ h_Q^{(\varepsilon)} :  Q\in \cS, \varepsilon \in \cA \}$
is a complete orthogonal system in $L^2(\bR^n) .$
Given 
$f \in L^p(\bR^n)$ 
\begin{equation}
\label{4maerz3}
 C_p^{-p} \|f \|_{L^p(\bR^n)}^p  \le \int\left
(\sum_{\varepsilon\in \cA,\,Q\in \cS}
 \la f,   h_Q^{(\varepsilon)} \ra^2  1_Q |Q|^{-2} \right)^{p/2} 
 \le 
C_p^p \|f \|_{L^p(\bR^n)}^p ,
\end{equation}
where $C_p \le Cp^2/(p-1). $
As is well known,  Wavelets, Calderon Zygmund operators and
Haar functions 
are related by $L^p$ equivalence. 
Any admissible  wavelet system 
$\{ \vp_Q^{(\varepsilon)}  : Q \in \cS, \varepsilon \in \cA \} $ is  
equivalent
to the Haar system $\{ h_Q^{(\varepsilon)} : Q \in \cS, \varepsilon
\in \cA \} $ in $ L^p(\bR^n) . $ 
For any choice of finite sums,
$$
f =\sum_{\varepsilon \in \cA,\, 
Q \in \cS}   a_Q^{(\varepsilon)} h_Q^{(\varepsilon)} 
\quad\text{and}\quad g =
\sum_{\varepsilon \in \cA , \,
Q \in \cS} a_Q^{(\varepsilon)} \vp_Q^{(\varepsilon)} , 
$$ 
we have equivalent norms
\begin{equation}
\label{5maerz1}
C(p ,\a,\d)^{-1} \left\| f
\right\|_{ L^p (\bR^n)} \le  
\left\| g
\right\|_ {L^p (\bR^n)} \le C(p ,\a,\d)
\left\|f
 \right \|_ {L^p (\bR^n)}. 
\end{equation} 
where  $C(p ,\a,\d) \le C(\a,\d)  p^2/(p-1) .$
See  \cite{coifme},
\cite{figsingular},  \cite{MR2157745} \nocite{zC87}.

\subsection{Rearrangements I} 
We present two quick applications of Semenov's theorem. The aim
is to estimate  series formed by  block bases of compactly supported
wavelets, Proposition \ref{19april011}.
\subsubsection*{Semenov's Theorem.}
Let $ \mu \in \bZ ^n . $  Semenov's
theorem \cite{MR2157745,pfxm97} asserts that the
rearrangement operator defined as the linear extension of 
$$ T_\mu :  h_Q^{(\varepsilon)} \to  h_{Q + \mu
  s(Q)}^{(\varepsilon)} $$
defines a bounded operator on $ L^p ( \bR^n ) $ with 
$$ \| T_\mu \|_p \le C_p \log (2 + |\mu| ) . $$ 
For our purposes the logarithmic dependence on $|\mu | $ 
is  crucial.

\begin{prop}\label{semenov}Let
$ \{ \vp_Q^{(\varepsilon)} :  Q\in \cS, \varepsilon \in \cA \}$
denote the wavelet  system defined by \eqref{structure011}. 
Let $\mu \in \bZ^n , $  and $f \in L^p ( \bR^n ) . $
Then
$$
\left\| ( \sum _{Q \in \cS } \la f , \vp_Q^{(\varepsilon)} \ra ^2 1 _{Q +
  \mu s(Q) } |Q|^{-2} ) ^{1/2} \right\|_p \le C(p,\a, \d) \log ( 2 + |\mu| ) \|f \| _p . $$  
\end{prop}

\proof
By Semenov's theorem and \eqref{4maerz3}, we have square function
estimates 
as follows
\begin{equation}\label{him10}
\left\| ( \sum _{Q \in \cS } \la f , \vp_Q^{(\varepsilon)} \ra ^2 1 _{Q +
  \mu s(Q) } |Q|^{-2} ) ^{1/2}\right \|_p \le  C(p) \log ( 2 + |\mu| )
\left\| ( \sum _{Q \in \cS } \la f , \vp_Q^{(\varepsilon)} \ra ^2 1 _{Q  }
|Q|^{-2} ) ^{1/2} \right\|_p.\end{equation}  
By \eqref{5maerz1} the Haar and wavelet systems are equivalent, so
that with \eqref{4maerz3}, the right hand side of \eqref{him10} is bounded by 
$$C(p,\a, \d) \log ( 2 + |\mu| ) \|f \| _p . $$
\endproof
\subsubsection*{Block bases of compactly supported Wavelets.} 
We consider again a   compactly supported  wavelet  system  in $L^2(\bR^n)$
$\{ \psi_K^{(\b)} : K \in \cS, \b \in \cA \} $
satisfying the structure conditions \eqref{13april_1}.

Let  $ \{c_K (Q) : K , Q \in \cS \} $ be a sequence of coefficients.
Define  
the following block basis, 
\begin{equation}\label{psitilde}
 \widetilde\psi_Q =   
\sum_{K\in  \cS_{j +k} }
      c_K(Q)   \psi_K^{(\b)} , \quad\quad Q \in \cS_j , \quad   j \in \bZ ,  \quad k \in \bN .
\end{equation} 
and form the operator, 
$$
S_0(f) = \sum_{Q \in \cS  }
\la f  , \vp_Q \ra  \widetilde\psi_Q    |Q|^{-1} . $$
Our aim is to prove that $S_0$ is bounded on $L^p(\bR^n)$
whenever  $ \{c_K (Q) : K , Q \in \cS \} $ satisfies \eqref{5nov011-1}.
To this end we split the block basis along integer translates of 
$Q$ and estimate with Semenov's theorem. To be precise,
 let  $ \mu \in \bZ^n . $ 
Then put 
 $$  \cA{ (Q, \mu )} = 
\{   K\in  \cS_{j + k}
 : K \sbe Q + \mu s(Q) \} \quad\text{ and }\quad \widetilde\psi_{Q, \mu } =   
\sum_{K\in   \cA {(Q, \mu ) } }
      c_K(Q)   \psi_K^{(\b)} . $$
Define next 
$$
S_\mu(f) = \sum_{Q \in \cS  } 
\la f  , \vp_Q \ra  \widetilde\psi_{Q, \mu }   |Q|^{-1} . $$

\begin{prop}\label{19april011}
Let $ \d > 0 , j \in \bZ, k\in \bN . $
Assume that the sequence of coefficients $\{c_K (Q) : K , Q \in \cS \} $
defining \eqref{psitilde} satisfy 
\begin{equation}\label{5nov011-1}
 |c_K (Q)| \le  \left( 1 + \frac{\dist( K , Q ) }{s(Q) }\right)^{
  -n(1+\d)} , \quad \quad Q  \in \cS_j   ,  K \in \cS _{j +k} . 
\end{equation}
Then for any  $ \mu \in \bZ^n ,$
$$ \|S_\mu \|_p \le C(p, \a, \d ) (1+|\mu |) ^{-n(1+\d)}  \log (2+|\mu |) , $$
and consequently
\begin{equation}\label{28ap011-B}
\|S_0 \|_p \le 
  C(p, \a,\d) . 
\end{equation}
\end{prop}
\proof 
Put
$$ \s^2(\widetilde\psi_{Q , \mu } ) := 
\sum_{K\in   \cA _{(Q, \mu ) } }
      |c_K(Q)|^2   |\psi_K^{(\b)}|^2 . $$
By \eqref{5nov011-1} 
 we have 
\begin{equation}\label{point011}
 \s^2(\widetilde\psi_{Q , \mu } )
\le 
C (1+|\mu |)^{-2n(1+\d)}
      \sum_{|\mu - \nu | \le C }1_{Q +
        \nu s(Q) } . 
\end{equation}
By the unconditionality of wavelet-bases  
\begin{equation}\label{28ap011-A} 
\left\|  S_\mu(f)\right\|_p
\le  C(p, \a,\d) \left\| (
\sum_{Q \in \cS  }
\la f  , \vp_Q \ra^2 \s^2( \widetilde\psi_{Q, \mu } ) |Q|^{-2}
)^{1/2}\right\|_p .
\end{equation}
By   \eqref{point011}
and Proposition~\ref{semenov}   the right hand side of 
\eqref{28ap011-A}  is
bounded by 
$$  C(p, \a,\d) (1 +|\mu |)^{-n(1+\d)}  \log( 2 + | \mu | ) \|f\|_p. $$
Since 
$$
S_0(f) = \sum_{\mu \in \bZ^n } S_\mu(f) $$
this gives  
$$ \|S_0 \|_p \le  C(p, \a,\d) \sum_{  \mu \in \bZ ^n}  \log( 2 + | \mu | ) (1 + |\mu |)^{-n(1+\d)}  $$
and  \eqref{28ap011-B} holds true.
\endproof

\subsection{Rearrangements II}\label{rearrangements}
We review here an auxiliary  result on a rearrangement 
operator $S$ that  induced by maping a dyadic cube to one of its dyadic predecessors.  
The operator was introduced and studied in  \cite{LMM}.
We  define
$S$ in \eqref{25jan0612}  and record its norm estimates. 
Let $\lambda \in \bN $ and let $Q\in \cS $ be a dyadic cube. 
The $\lambda -th $ dyadic predecessor of $Q,$ denoted 
$Q^{(\lambda)}, $ is given by the relation 
$$Q^{(\lambda)} \in \cS , \quad |Q^{(\lambda)}| = 2^{n\lambda}|Q|,\quad 
Q \sb Q^{(\lambda)}. $$
Let $\tau :\cS \to \cS $ be the map that associates to 
each $ Q \in \cS $ its $\lambda -th $ dyadic predecessor. 
Thus
$$\tau(Q) =  Q^{(\lambda)},\quad Q \in  \cS .$$
Clearly $\tau  :\cS \to \cS $ is not injective.
We canonically split $ \cS = \cQ_1 \cup \cdots \cup  \cQ_{2^{n \lambda}} $
such that the restriction of $\tau $ to  each of the collections
 $ \cQ_k ,$
is injective:
Given $ Q \in \cS ,$ form
$$ \cU(Q)  = \left\{ W \in \cS : W^{(\lambda)} = Q \right\} . $$
Thus $ \cU(Q) $ is a covering of $Q$ and 
contains exactly  $2^{n \lambda}$ pairwise disjoint dyadic cubes.
We enumerate them,
rather arbitrarily, as
$ W_1(Q) , \dots ,  W_{2^{n \lambda} }(Q). $
For $ 1 \le k \le  2^{n \lambda} ,$  define 
$$ \cQ_k = \left\{ W_k(Q) : Q \in \cS \right\} . $$
Note that  $\tau :\cQ_k \to \cS $ is a bijection, and 
$$ \tau (  W_k(Q) ) = Q , \quad\quad  W_k(Q)\in \cQ_k,
\quad Q \in \cS  . $$

Let  $ 1 \le k \le  2^{n \lambda} .$ Let 
$\{ F_Q^{(k)} : Q \in \cS \} $ be any family of functions satisfying
$ \int F^{(k)}_Q (x) dx = 0 $ and the following 
structural conditions: There exists $C>0 $ $ \d > 0 $ and $ 0 < \a \le
1 $ so that for each $Q \in \cS $
\begin{subequations}\label{11juli12}
\begin{equation}\label{11juli12a}
\quad\quad  | F^{(k)}_Q (x)| \le C\left( 1 + \frac{\dist(x, Q) }{s(Q) }
 \right)^{-n(1+\d)},  \end{equation}  
and for $  |x-t | \le  s(Q) ,$
\begin{equation}\label{11juli12b}
 | F^{(k)}_Q (x) - F^{(k)}_Q (t)| 
\le  C s(Q)^{-\a}|x-t|^\a \left( 1 + \frac{\dist(x, Q) }{s(Q) }
\right)^{-n(1+\d)} . 
\end{equation}
\end{subequations}
We emphasize that $  F_Q^{(k)}$  may  depend
on $k ,$ by contrast 
  the structural conditions \eqref{11juli12} 
 are independent of  the value of
$k . $ Define the operator $S$ by the equation
\begin{equation}
\label{25jan0612}
 S(g) = \sum_{k = 1 }^{ 2^{n \lambda} }
\sum_{Q \in  \cQ_k } \left\la g ,  F_{\tau(Q)}^{(k)}\right\ra \vp_Q |Q|^{-1},
\end{equation}
where $\{\vp_Q\}$ is an admissible  wavelet system satisfying \eqref{structure011}.
The operator $S$ is the transposition of the rearrangement operator 
defined by $\tau $ followed by a  Calderon Zygmund Integral.
The next theorem records the 
 operator norm of $ S,$ particularly  its joint $(n,\lambda) -$dependence, 
on  $L^p (\bR^n).$ 
\begin{theor}
\label{11julitheorem1}
Let $1 < p < \infty . $
The operator  $ S$ defined by \eqref{25jan0612} is bounded  $L^p (\bR^n)$
The 
norm estimates depend on the value of $\lambda \in \bN$ and the dimension of the 
ambient space $\bR^n $ as follows:
\begin{equation}\label{11juli20}
 \|S\|_p \le   C(p, \a, \d)\lambda ^{1/2}   2^{n \lambda} .
\end{equation}
\end{theor}
\proof Just transfer Theorem 5.2 in \cite{LMM}
from compactly supported wavelets to those satisfying \eqref{11juli12}.
\endproof

\section{Proof of Theorem~\ref{th2a}.}
\label{proof}
In this section we prove Theorem~\ref{th2a}.
The sub-section ~\ref{sub1} is devoted to the estimates \eqref{8jan3} for the operator
$T_\ell^{(\varepsilon)},$ $\ell\ge 0.$ 
Thereafter we discuss the reduction of the 
estimates \eqref{81jan3} for $T_\ell^{(\varepsilon)} R^{-1}_{i_0},$ $ \varepsilon 
\in \cA_{i_0},$ to those of 
$T_\ell^{(\varepsilon)}.$

\subsection{Estimates for $T_\ell^{(\varepsilon)} . $}

We prove here \eqref{8jan3} asserting
that $T_\ell^{(\varepsilon)} ,\,\ell \ge 0 $ satisfies the norm estimates
\begin{equation}
\label{28mai3}
\| T_{\ell}^{(\varepsilon)} \|_p \le   C(p, \a, \d)2^{-\ell\a } . 
\end{equation}
We do this by performing a further decomposition of the operator
$T_\ell^{(\varepsilon)}. $ 
\subsubsection*{The    Decomposition of  $T_\ell^{(\varepsilon)},
  \,\ell \ge 0 .$}

We 
 decompose the 
operator $T_\ell^{(\varepsilon)},$ $\ell \ge 0 $ into a series of operators 
$T_{\ell, m}, m \in \bZ$  using 
 a compactly supported 
wavelet system 
$\{ \psi_K^{(\b)} : K \in \cS, \b \in \cA \} .$
We assume that 
 $\{ \psi_K^{(\b)} /\sqrt{| K|}
\} $
is an orthonormal basis in $L^2(\bR^n),$
satisfying  $\int \psi_K^{(\b)} = 0 $ and
the structure conditions,
$$
 \supp\psi_K^{(\b)} \sbe C\cdot K,
\quad \quad |\psi_K^{(\b)}| \le C ,
\quad\quad \Lip(\psi_K^{(\b)}) \le C \diam(K)^{-1}.
$$
 We suppress the superindeces  $(\b)$   and,
in place  of  $\{\psi_K^{(\b)}\}$ we write 
just $\{ \psi_K \}.$ 

 Fix $m\in\bZ , j \in \bZ ,$ and $  Q \in\cS_j .$ Put 
\begin{equation}\label{blockbasis}
\widetilde \psi_{Q}=
\sum_{K\in \cS_{j +\ell+m}}
 \left \la \Delta_{j+\ell}(\vp^{(\varepsilon)}_{Q}),\psi_K
\right\ra \psi_K |K|^{-1} 
\end{equation}
and 
\begin{equation}
T_{\ell,m}(f)= 
\sum_{Q \in \cS }
\left\la f, \widetilde\psi_{Q  }
\right \ra
\vp^{(\varepsilon)}_{Q}
|Q|^{-1} .
\end{equation} 
By construction, $$T_\ell = \sum_{m = -\infty}^\infty T_{\ell,m} . $$

For fixed $ \ell \ge 1 $ we consider below  three cases
$$-\infty < m \le -\ell - 1 ,\quad
 -\ell \le m \le 0 ,\quad\text{ and }\quad m \ge 0 . $$
We  prove accordingly that 
\begin{equation}\label{dez511a}
\sum^{-\ell - 1}_{m=-\infty}||T_{\ell,m}||_p  \le
   C(p, \a, \d) 2^{-\ell \a},
\quad \quad 
\sum_{m=-\ell}^0
\| T_{\ell,m} \|_p \le   C(p, \a, \d) 2^{-\ell \a} \end{equation}
and 
\begin{equation}\label{dez511b}
\sum_{m= 0}^\infty
\| T_{\ell,m} \|_p \le   C(p, \a, \d) 2^{-\ell \a}. 
\end{equation}
The estimates  \eqref{dez511a} and \eqref{dez511b} 
yield  $\| T_{\ell} \|_p \le   C(p, \a, \d)
2^{-\ell \a}  $ as claimed.

\label{sub1}

\begin{prop} \label{pr1}
Let $1 < p < \infty .$ Let  $ \d > 0 $ and  $\a > 0 $ be  fixed  in
the definition of the 
admissible wavelet system.  Put  $ \gamma_0 = \min \{n\d /2, 1 \} . $  
For   $\ell \ge 0$,  and   $m < -\ell$ the operator $T_{\ell,m} $
satisfies the norm estimate
\begin{equation}
\label{23maerz1}
\| T_{\ell , m } \|_p\le  C(p, \a, \d) 2^{-\ell \a} 2^{-|m+ \ell |
  \gamma_0}\sqrt{|m +\ell|} .
\end{equation}
\end{prop}
\proof 
Let $j \in \bZ$ and fix
 a dyadic cube $Q \in \cS_j . $
Since $ \ell + m  < 0 $  there exists a unique cube  $K_0 \in \cS_{j
  +\ell+m} ,$  so that $ Q \sbe K_0 .$ 
Lemma \ref{12april_3} gives the point-wise estimates 
\begin{subequations}\label{27ap011B} 
\begin{equation}\label{27ap011C} 
|\widetilde \psi_{Q} (t)| \le C  2^{-\a\ell} 2^{-|m+\ell | (n + \g_o) }
\left(1 +\frac{\dist(t , K_0)}{s(K_0)} \right)^{-n(1+\d)}, 
\end{equation}
\begin{equation}\label{27ap011D} 
|\nabla \widetilde \psi_{Q} (t)| \le C \diam (K_0) ^{-1} 2^{-\a\ell}
2^{-|m+\ell | (n + \g_o) }
\left(1 +\frac{\dist(t , K_0)}{s(K_0)} \right)^{-n(1+\d)}. 
\end{equation}
\end{subequations}

We next  invoke rearrangement operators. 
Let $\tau : \cS \to \cS $ be the map that associates 
to $ Q \in \cS $ its $|m +\ell | -th $ dyadic predecessor,
denoted $Q^{|m +\ell |} . $
Thus 
$$ \tau( Q ) = Q^{|m +\ell |} . $$
In sub-section ~\ref{rearrangements} we defined the canonical splitting of $\cS$
as 
$$\cS = \cQ_1 \cup \dots \cup \cQ_{2^{n|m +\ell |}},$$
so that for each fixed  $ k \le 2^{n|m +\ell |},$
the map
$\tau : \cQ_k \to \cS $ is a bijection.
Fix now $ k \le 2^{n|m +\ell |}$ and define the family of functions
$\{ F_W ^{(k)} : W \in \cS \} $ by the equations
 
\begin{equation}\label{7dez11d}
  F_{\tau(Q)} ^{(k)} =   2^{\a\ell} 2^{|m+\ell| (n+\g_0) }
  \widetilde \psi_{Q}, \quad\quad Q \in  \cQ_k .
\end{equation}
Let  $A = 2^{n|m +\ell  |} $ and  define the rearrangement operator $ S $ by 
$$
 S(f) = \sum_{k = 1}^A \sum_{Q \in \cQ_k }
\left\la f,    F_{\tau(Q)} ^{(k)}    \right\ra \vp_{Q}|Q|^{-1}.
$$

By \eqref{27ap011B},      $\{ F_W ^{(k)} : W \in \cS \} $ satisfies the structure 
estimates \eqref{11juli12}. Apply  Theorem~\ref{11julitheorem1} 
with  $\lambda = |m +\ell| , $
 to obtain
$$ \|S\|_p \le
 C(p, \a, \d) 2^{n |m +\ell | }  \sqrt{|m+\ell|}   
$$
Hence  with \eqref{7dez11d} we get
\begin{equation}\label{11juli21}
 \|T_{\ell, m }(f)\|_p \le  C(p, \a, \d)  2^{-\a\ell}  2^{-|m+\ell| (n+\g_0) } 
\| S(f)\|_p \le    C(p, \a, \d)  2^{-\a\ell}  2^{-|m+\ell| \g_0 } \sqrt{|m+\ell|} 
\|f\|_p  . 
\end{equation}
\endproof

We treat next the case $m \ge  0$ and   $\ell\ge 0.$ 
Here we estimate
the  transposed
operator of $T_{\ell,m} $
which  is given by
\begin{equation}\label{10dez111}
T_{\ell,m}^*(f)= 
\sum_{Q \in \cS }
\left\la f, \vp^{(\varepsilon)}_{Q} \right \ra \widetilde\psi_{Q  }
|Q|^{-1}.
\end{equation} 

\begin{prop} \label{pr3}
Let  $1 < p < \infty.$ 
For   $m\ge 0$ and   $\ell\ge 0.$
Then 
\begin{equation}
\label{18feb20}
\| T_{\ell,m} \|_p \le  C(p, \a, \d)2^{-m} 2^{-\ell \a}
\end{equation}
\end{prop}
\proof
Fix $\ell \ge 0$   
and  $   m\ge 0.$
Let $j \in \bZ $ and choose a 
 dyadic cube $Q \in \cS_j . $
The structure estimates for $\Delta_{j+\ell}(\vp^{(\varepsilon)}_{Q})$
in Lemma~\ref{basiclemma}
translate into coefficient estimates as follows.
If $K \in \cS_{j + \ell + m } ,$ then 
\begin{equation}
\label{18feb3}  
\begin{aligned}
|\la\Delta_{j+\ell}(\vp^{(\varepsilon)}_{Q}), \psi_K\ra|\cdot|K|^{-1}
&\le C 2^{-m} 2^{-\a \ell} \left( 1 + \frac{\dist(K, Q )}{ s(Q)}\right) ^{ -n(1+\d)}. 
\end{aligned}
\end{equation}
Using \eqref{18feb3} and applying  Proposition~\ref{19april011} to  $T_{\ell,m}^*$
gives the norm estimate
$$
\| T_{\ell,m} \|_p \le  C(p, \a, \d)2^{-m} 2^{-\ell \a} . $$
\endproof

Next consider $\ell \ge 0,$ $-\ell\le m\le 0.$
The ingredients  of the previous proof are applied again.
\begin{prop}\label{pr2}
Let  $1 < p < \infty.$ 
Let    $\ell \ge 0$ and  $-\ell\le m\le 0$ then  
\begin{equation}
\label{23maerz10}
\| T_{\ell,m} \|_p \le   C(p, \a, \d) 2^m 2^{-\ell \a} 
\end{equation}

\end{prop}
\proof
We estimate the transposed operator
$T_{\ell,m}^*$ given by \eqref{10dez111}.
Fix $\ell \ge 0$ and  $-\ell\le m\le 0.$
Let $j \in \bZ $ and choose 
 dyadic cubes $Q \in \cS_j  $  and 
  $K \in \cS_{j + \ell + m } .$
Next we apply  Lemma~\ref{scalar-prod} and get 
$$
|\la\Delta_{j+\ell}(\vp^{(\varepsilon)}_{Q}),
\psi_K\ra|\cdot|K|^{-1}
\le C 2^m 2^{-\a \ell} \left( 1 + \frac{\dist(K, Q )}{ s(Q)}\right) ^{
  -n(1+\d)}.
$$
Applying 
Proposition~\ref{19april011}  to $T_{\ell,m}^*$ gives
$\| T_{\ell,m} ^*\|_p \le   C(p, \a, \d) 2^m 2^{-\ell \a}.$
\endproof

\subsubsection*{The proof of Theorem \ref{th2a}. Part 1.} 
The estimate
\eqref{8jan3} 
is now obtained as follows:
The assertions of Proposition \ref{pr1}, Proposition \ref{pr3} and Proposition \ref{pr2}  
imply
\begin{equation}
\label{28feb1}
\sum^{\infty}_{m=-\infty}||T_{\ell,m}||_p  \le
   C(p, \a, \d) 2^{-\ell \a},
\end{equation}
Since \begin{equation}\label{8jan9}
T_{\ell}^{(\varepsilon)} (f) = \sum^\infty_{m=-\infty}T_{\ell,m} (f),
\end{equation}
we get
\eqref{8jan3}.
\endproof
\label{sub3}

\subsection{Estimates for $T_{\ell}^{(\varepsilon)}
R^{-1}_{i_0}.$ }
\label{sub4}
Here we prove \eqref{81jan3}.
We fix $\ell \ge 0,$  $1\le i_0 \le n ,$ and  
$ \varepsilon \in \cA_{i_0}.$ 
We now prove the  obtain  norm estimates for
$T_{\ell}^{(\varepsilon)}R^{-1}_{i_0} $
 {\em by reduction} to the estimates 
for
the operator $T_{\ell} ^{(\varepsilon)} $

Let $ j \in \bZ  $ and $ Q \in \cS_j .$ 
Recall that we put
$$
k_Q^{( \ell, i )} = 
\Delta_{j+\ell}\left( \bE_{{i_0}}\pa{_i}\vp_Q^{( \varepsilon )}\right).
$$
\begin{prop}
\label{12juli4}
Let $1 < p < \infty. $ Let $ 1 \le i \ne i_0 \le n $ and $\varepsilon \in
\cA_{i_0}.$
For $\ell \ge 0 $ the operator $X$ defined by  
$$X(f)=
\sum_{Q \in \cS}
\la
 f,   k_Q^{( \ell, i )}
\ra \vp^{(\varepsilon)}_Q |Q|^{-1}, $$
satisfies the norm estimates
\begin{equation}
\label{23maerz5}
||X||_{p}\le 
 C(p, \a, \d)2^{+\ell -\a\ell }  
\end{equation}
\end{prop} 

\proof
It remains to compare the structure conditions  \eqref{27mai3}     for the system 
$ k_Q^{( \ell, i )}$ defining $X$ with those for $ f_Q^{( \ell, i )}$ defining  the operator 
$T_\ell^{(\varepsilon)} $
This gives 
$$ 
||X||_{p}\le  2^{\ell}\|T_\ell^{(\varepsilon)} \|_p ,
$$
which implies  \eqref{23maerz5}.
\endproof

\subsubsection*{The Proof of Theorem \ref{th2a} Part 2.}

The  estimate
\eqref{81jan3} is obtained as follows: 
Proposition~\ref{12juli4} in combination with the norm estimate \eqref{8jan3}
and the representation\eqref{21maerz1}
imply  that for $\ell > 0 ,$
$$||T_\ell^{(\varepsilon)} R^{-1}_{i_0}||_{p}\le  C(p, \a, \d)2^{+\ell -\a\ell }
 .
$$

\section{Proof of Theorem~\ref{th2b}.}
\label{simple}
In this section we prove Theorem~\ref{th2b}.
For  $\ell \le 0$ we obtain the norm estimates for 
$T_\ell ^{(\varepsilon)}R^{-1}_{i_0}$ and $T_\ell^{(\varepsilon)} $
by the same method.
 Let 
$ i \ne i_0 $ and  $ \varepsilon \in \cA_{i_0}$ and let $\ell \le 0 . $

\subsubsection*{The Proof of Theorem  \ref{th2b}. }
 Let $ 1 < p < \infty $ , $\ell \le 0  $ and $ \varepsilon \in
 \cA_{i_0} .$
We show that then 
$$
||T_\ell ^{(\varepsilon)}||_{p} + ||T_\ell^{(\varepsilon)} R^{-1}_{i_0}||_{p} 
\le 
C(p, \a ,\d) 2^{-|\ell|} |\ell|
$$
\nocite{MR2157745}
We use the representations \eqref{27mai7} and \eqref{21maerz1}, and recall 
that we put  
$$
  f^{(\varepsilon)}_{Q,\ell} = \Delta _{j+\ell} ( \vp_Q^{(\varepsilon)}) ,
\quad\text{and}\quad k_Q^{( \ell, i )} = 
\Delta_{j+\ell}( \bE_{{i_0}}\pa{_i}\vp_Q^{( \varepsilon )}), 
\quad Q \in \cS_j.
$$
We showed in  Lemma \ref{9dez11a} that  $ \{  k_Q^{( \ell, i )} : \,Q \in \cS \, ,\ell
\le 0 \} $
satisfies  conditions \eqref{27mai4}. 
 It is easy to see that 
also  the 
family $ \{  f^{(\varepsilon)}_{Q,\ell}: \,Q \in \cS \, ,\ell \le 0 \} $ 
satisfies 
 the  structural conditions \eqref{27mai4}.

Now  we
 choose $\{g_{Q,\ell}: Q \in \cS\} $ satisfying the structure conditions 
\eqref{27mai4}. Define the operator
$$
X(f) = \sum_{Q \in \cS } \la u ,g_{Q,\ell}\ra 
 \vp_Q^{(\varepsilon)} |Q|^{-1}. $$
In view of the preceding discussion, the $L^p $ estimates for $X$
will 
apply to both 
$T_\ell ^{(\varepsilon)}$ and $T_\ell^{(\varepsilon)} R^{-1}_{i_0}.$

To estimate $X ,$ we consider again the rearrangement 
$\tau :\cS \to \cS $ that maps $ Q\in \cS $ 
to its $|\ell| -th $ dyadic predecessor.
Let $ \cQ_1 , \dots ,\cQ_{2^{n|\ell|}}$ be the canonical splitting
of $ \cS$ so that for fixed 
$ k \le   2^{n|\ell|} $ the map 
$\tau :\cQ_k \to \cS $
is bijective. Fix $ k \le   2^{n|\ell|} .$
Determine the family 
$\{ F_W ^{(k)} : W \in \cS \} $ by the equations
\begin{equation}\label{28jan062} 
  F_{\tau(Q)} ^{(k)} = 2^{(n+1) |\ell |} g_{Q,\ell}, \quad\quad Q \in  \cQ_k . \end{equation}
Define the operator
$$
S(u) = \sum_{k = 1 } ^{2^{n |\ell |}} \sum_{Q \in \cQ_k } 
\left\la u , F_{\tau(Q)} ^{(k)}\right\ra
 \vp_Q^{(\varepsilon)} |Q|^{-1} .
$$
Apply  Theorem~\ref{11julitheorem1} to $S$ with $\lambda = | \ell | . $
This yields 
\begin{equation}
\label{25jan0617}
 \| S\|_p \le  C(p, \a ,\d) 2^{n |\ell|}|\ell|
\end{equation}
Comparing the structure conditions gives  
\begin{equation}
\label{18feb31}
 \|X \|_p  \le 
C(p, \a ,\d) 2^{-(n+1)|\ell|}
\|S\|_p .
\end{equation}
Consequently, our  upper bounds for  $\|X \|_p $
follow  from \eqref{25jan0617}.
Indeed,
$$\|X \|_p 
\le     C(p, \a ,\d) 2^{-|\ell|}  |\ell|
$$

The above estimate for $X$ gives finally 
$$
||T_\ell ^{(\varepsilon)}||_{p} +
||T_\ell ^{(\varepsilon)}R^{-1}_{i_0}||_{p} 
\le      C(p, \a ,\d)2^{-|\ell|}  |\ell| .$$
\endproof

\nocite{jpw1}
 \bibliographystyle{abbrv}
\bibliography{compensated}

\noindent
{\bf  AMS Subject classification: 49J45 , 42C15, 35B35}\\

\noindent
{\bf  Addresses:}\\

\noindent
 Paul F.X. M\"uller\\
Institut f\"ur Analysis \\
J. Kepler Universit\"at\\
 A-4040 Linz \\
E-mail: paul.mueller@jku.at
\vskip 0.3cm
\noindent 
Stefan M\"uller\\
Universit\"at Bonn\\
Endenicher Allee 60\\
D 53115 Bonn Germany\\
E-mail: stefan.mueller@hcm.uni-bonn.de
\vskip 0.5cm
\noindent
\end{document}